\documentclass[preprint,number, draft]{elsarticle}
\usepackage{lineno}
\usepackage{amsfonts,amssymb,amsmath}

\newtheorem{thm}{Theorem}[section]
\newtheorem{lem}{Lemma}[section]
\newtheorem{prop}{Proposition}[section]
\newtheorem{cor}{Corollary}[section]
\newtheorem{ex}{Example}[section]
\newdefinition{newdef}{Definition}[section]
\newdefinition{rmk}{Remark}[section]
\newproof{pf}{Proof}

\modulolinenumbers[5]
\journal{Linear Algebra and its Applications}

\begin{document}

\begin{frontmatter}

\title{On Wiener -- Hopf factorization of scalar polynomials}

\author{Victor M. Adukov}
\address{Faculty of Mathematics, Mechanics and Computer Science,\\ Institute of Natural Sciences,\\South Ural State University
(national research university), Chelyabinsk, Russia}
\ead{adukovvm@susu.ac.ru}

\begin{abstract}
In the work we propose  an algorithm for a Wiener -- Hopf factorization of scalar polynomials
based on notions of indices and essential polynomials. The algorithm uses computations with
finite Toeplitz matrices and permits to obtain coefficients of both factorization factors simultaneously. Computation aspects of the algorithm are considered. An a priory estimate for the condition number of the used Toeplitz matrices  is obtained.
Upper bounds for the accuracy of the factorization factors are established. 
All estimates are effective. 
\end{abstract}

\begin{keyword}
{ Wiener -- Hopf factorization \sep polynomial factorization \sep Toeplitz matrices}
\MSC[2010]47A68 \sep 12D05 \sep 12Y05
\end{keyword}

\end{frontmatter}


\section{Introduction}

The Wiener -- Hopf technique is a powerful method used in various areas of mathematics, mechanics and mathematical physics (see, e.g., \cite{DZ}, \cite{GF}, \cite{LA}).
The core of the method is the Wiener -- Hopf factorization problem for matrix functions
(or the Riemann boundary-value problem) \cite{CG,GF}. In the scalar case there is an explicit formula for a solution of the problem \cite{Gakhov}, but its use for numerical computations is quite difficult. Therefore to factorize a scalar function, it is usually approximated by a rational function (see, e.g., \cite{Kisil}). Thus the scalar factorization reduces to the polynomial factorization. 

The matrix case is more complicated because there are no explicit formulas for the factorization factors and for important integer invariants of the problem such as partial indices. Moreover, it is difficult to develop approximate methods of the factorization since
the matrix factorization is unstable. For this reason it is very important to find cases when the problem can be effectively or explicitly solved. 
Almost all cases of constructive factorization known to date are considered in the review
\cite{RM}. In~\cite{Adukov92} an explicit method for the factorization of meromorphic matrix functions was proposed. In this paper 
the factorization problem is called
explicitly solved if it is reduced to the factorization of given scalar functions and to solving of finitely many finite systems of linear algebraic
equations. In~\cite{Adukov08a} some classes of matrix functions are listed for which the factorization problem can be explicitly reduced to the factorization of analytic matrix functions.

In this paper, we start a project to develop algorithms of the factorization and their implementations for these classes of matrix functions.  We propose to consider the factorization of matrix polynomials, analytic and meromorphic matrix functions, triangular matrix functions and matrix functions with  one non-meromorphic  
row.   Some preliminary results in  this direction were 
obtained in~\cite{Adukov08}.

The first stage in solving of the factorization problem for the above mentioned classes is the factorization of scalar functions that can be reduced to the polynomial factorization. 

We consider coefficients of a polynomial $p(z)$ as initial data and coefficients of its factors as output data of the problem. 
The naive method of the polynomial factorization is the following way. First we compute roots of a polynomial $p(z)$ and perform their separation. Second we find the coefficients of the factors of $p(z)$ by their roots. 
It is well known that the roots of a polynomial are in general not well-condition functions of its coefficients (see, e.g., \cite{Gau}), and 
coefficients of a polynomial are also 
not well-condition functions of its roots~\cite{Uh}.
 
 The latter means that, in general, we can not solve numerically the polynomial factorization problem by the naive way.

 Nevertheless there exist numerical methods for solving of this problem. The basic works in this direction are cited in~\cite{BB}.
It should be especially pointed out the articles~\cite{BB}--\cite{BGM2}, where 
finding of 
the polynomial factorization was based on computations with Toeplitz matrices.

Our approach resembles the method  proposed by D.A. Bini and A. B\"ott\-cher~\cite{BB} in 
{\bf algorithm 3}.
As the authors we solve systems with finite Toeplitz matrices consisting of Laurent coefficients of the function $1/p(z)$. 
The main difference between the method of D.A. Bini and A. B\"ottcher and ours is that 
{\bf algorithm~3} permits to obtain only coefficients of the factor $p_-(z)$ for $p(z)$, whereas we find coefficients of both factors simultaneously. This is important because the polynomial division is, in general, not well-condition operation in numerical computations.
Moreover, we use a different technique that can be extended to the factorization problem for analytic functions. In this case
we can obtain all coefficients of the polynomial factor and a required number of Taylor coefficients of the analytic factor.

The paper is organized as follows. In section 2 we consider the setting of the polynomial factorization problem and formulate some results on used norms of polynomials and their estimates. Section 3 contains basic tools for solving of the problem. Here we introduce 
notions of indices and essential polynomials in terms of which the problem will be solved.
In Section 4 we prove the basic results on the factorization of a polynomial $p(z)$. In Section 5 some computation aspects of the factorization problem are considered. Here we obtain an a priory estimate for the condition number of the used Toeplitz matrices and establish 
upper bounds for the accuracy of the factorization factors. 
All estimates are effective.
Section 6 contains an algorithm and numerical examples.
\section{Preliminaries}

Let $p(z)=p_0+p_1 z+\cdots+p_{\nu }z^{\nu}$ be a complex polynomial of  degree $\nu>1$ and $p_0\ne 0,  p_\nu=1$.
We suppose that $p(z)\ne 0$ on the unit circle $\mathbb T$, hence $p(z)\ne 0$ on a  closed circular annulus
$K:=\left\{ z\in \mathbb C : r\leq |z|\leq R\right\}$ for some $0<r<1<R<\infty$. Denote by 
$\varkappa = {\rm ind}_{\mathbb T}\, p(z)$ the index of $p(z)$ with respect to  $\mathbb T$,
i.e. the number of zeros of the polynomial inside the unit circle.  Let $\xi_j$, $j=1,\ldots,\nu$, be the zeros of
$p(z)$ and 
$$
0<|\xi_1|\leq\ldots\leq|\xi_\varkappa|<r<1<R<|\xi_{\varkappa+1}|\leq\ldots\leq|\xi_\nu|.
$$
Here the zeros are counted according to their multiplicity.

In the work we will consider the following factorization of $p(z)$
\begin{equation}\label{factp0}
p(z)=p_1(z) p_2(z),
\end{equation}
where $
p_1(z):={(z-\xi_1)\cdots(z-\xi_\varkappa)}, \ p_2(z):=(z-\xi_{\varkappa+1})\cdots(z-\xi_\nu).$
Denote $p_-(z)=\frac{p_1(z)}{z^\varkappa}$, $p_+(z)=p_2(z)$.
Then the representation

\begin{equation}\label{factWp}
p(z)=p_-(z)z^\varkappa p_+(z), \  \ |z|=1,
\end{equation}
is the Wiener -- Hopf factorization of $p(z)$  normalized by the condition $p_-(\infty)=1$.

Throughout this paper, $\|x\|$ means the
 H\"older $1$-norm
$$
\|x\|=|x_1|+\cdots+|x_k|,
$$
$x=(x_1,\ldots,x_k)^T\in\mathbb C^k$.
 A norm $\|A\|$ of a matrix 
$A\in{\mathbb C}^{\,\ell\times k}$
is always the induced norm
$$
\|A\| = \max_{1\le j\le k}\sum_{i=1}^\ell|A_{ij}|.
$$

Respectively, the norm of a polynomial $p(z)=p_0+p_1 z+\cdots+p_{\nu} z^\nu$ is the norm of the vector $(p_0,p_1,\ldots,p_\nu)^T$.
For $p(z)$ we will also use the maximum norm
$$
\|p\|_C=\max_{z\in \mathbb T}|p(z)|
$$
on the unit circle $\mathbb T$. 

The norms $\|\cdot\|$ and $\|\cdot\|_C$ are equivalent.
Clearly, $\|p\|_C\leq\|p\|$. Since for $p(z)$ it is fulfilled the equality
$$
\sum_{k=0}^\nu|p_k|^2=\frac{1}{2\pi}\int_0^{2\pi} |p(e^{i \varphi})|^2d\varphi, 
$$
we have $\|p\|\leq \sqrt{\nu+1}\,\|p\|_2\leq\sqrt{\nu+1}\,\|p\|_C.$ Thus,
\begin{equation}\label{equnorm}
\|p\|_C\leq\|p\| \leq\sqrt{\nu+1}\,\|p\|_C.
\end{equation}

In order to study stability of the factorization problem, we will need estimates for the norm of  inverses of some Toeplitz matrices. Such estimates will be obtain in terms of 
$\|p_1\|\,\|p_2\|$, where $p_1(z)$, $p_2(z)$ are the factorization factors of $p(z)$. 
To get effective estimates, it will be required to estimate $\|p_1\|\,\|p_2\|$ via $\|p\|$.

Let $q(z)=q_1(z)q_2(z)$, where $q_1(z)$, $q_2(z)$ are arbitrary monic complex polynomials and $\nu =\deg q$. It is obvious that
$$
\|q\|\leq \|q_1\|\,\|q_2\|.
$$
In the work \cite{Boyd}, D.W. Boyd proved the following inequality
$$
\|q_1\|_C\,\|q_2\|_C\leq \delta^\nu \|q\|_C,
$$
where $\delta =e^{2G/\pi}=1.7916228120695934247\ldots$ and $G$ is Catalan's constant. The inequality is asymptoticaly sharp as $\nu\to \infty$.

Taking into account (\ref{equnorm}), in our case we obtain
\begin{equation}\label{estmulti}
\|p\|\leq\|p_1\|\,\|p_2\|\leq \delta^\nu\sqrt{(\varkappa+1)(\nu-\varkappa+1)}\,\|p \|.
\end{equation}

However, the exponential factor $\delta^\nu$ can overestimate the upper bound.
In some special cases we can obtain more precise estimates. For example, 
this can be done for a so-called spectral factorization of polynomials.

\begin{prop}\label{spcase}
Let $p(z)=\sum_{j=0}^{2m}p_jz^j$ be a complex polynomial of degree $2m$ such that $ p_{2m-j}=\bar p_j$ for $j=0,\ldots,m$, $p_0=1$. Suppose that $p(z)\ne 0$ on $\mathbb T$.

If $p(z)=p_1(z)p_2(z)$ is factorization~(\ref{factp0}), then
$$
\|p_1\|\,\|p_2\|\leq(m+1)\|p\|.
$$
\end{prop}
\begin{pf}
It is clear that if $\xi_k$ is a root of $p(z)$, then so is $1/\bar \xi_k$. 
Hence $m$ is the index of $p(z)$.

Let $|\xi_k|<1$ for $k=1,\ldots,m$ and
$$
p_1(z)=(z-\xi_1)\cdots(z-\xi_m), \ p_2(z)=(z-1/\bar\xi_{1})\cdots(z-1/\bar\xi_m)
$$
are the factorization factors of $p(z)$. On the unit circle we have
$$
\overline{p_1(t)}= (-1)^m(\bar\xi_1\cdots\,\bar\xi_m)t^{-m}p_2(t)=
\overline{p_1(0)}t^{-m}p_2(t), \ |t|=1.
$$
It follows from this that 
$\|p_1\|_C=|p_1(0)|\,\|p_2\|_C$ and
$$
p(t)=\frac{t^m}{\overline{p_1(0)}}\,p_1(t)\overline{p_1(t)}.
$$

Thus we have
$$
\|p\|_C=\frac{1}{|p_1(0)|}\|p_1\|_C^2=\|p_1\|_C\|p_2\|_C.
$$ 
Applying inequality (\ref{equnorm}), we arrive at the assertion. $\ \ \Box$
\end{pf}

It is possible that the equality $\|p\|=\|p_1\|\,\|p_2\|$ holds.

\begin{prop}\label{spcase1}
Let $p(z)=\sum_{j=0}^{2m}p_jz^j$ be a real polynomial of degree $2m$ and $p_{2m-j}=p_j$ for $j=0,\ldots,m$, $p_{0}=1$. Suppose that $p(z)\ne 0$ on $\mathbb T$ and all roots of $p(z)$ have negative
real parts.

If $p(z)=p_1(z)p_2(z)$ is factorization~(\ref{factp0}), then
$$
\|p\|=\|p_1\|\,\|p_2\|.
$$
\end{prop}
\begin{pf}
 In this case, if $\xi_k$ is a root of $p(z)$, then so are $\bar\xi_k$
and $1/\xi_k$. 

Let $|\xi_k|<1$ for $k=1,\ldots,m$. Then 
$$
p_1(z)=(z-\xi_1)\cdots(z-\xi_m), \ p_2(z)=(z-1/\xi_{1})\cdots(z-1/\xi_m).
$$

If $z=\xi$ is a complex root of $p(z)$, then $(z-\xi)(z-\bar\xi)$ has  positive  coefficients. Obviously, $z-\xi$ has also positive  coefficients for a real root $\xi$.
Then $p_1(z)$ has positive coefficients as a product of such polynomials.
Similar statement holds for $p_2(z)$ and $p(z)$.

Thus we have
$\|p\|=p(1)=p_1(1)p_2(1)=\|p_1\|\,\|p_2\|$. $\ \ \Box$
\end{pf}

Since for a given polynomial $p(z)$ more precise estimates can exist,  we will used the inequalities
\begin{equation}\label{estmulti1}
\|p\|\leq\|p_1\|\,\|p_2\|\leq \delta_0\|p \|,\ \ 1\leq\delta_0\leq\delta^\nu\sqrt{(\varkappa+1)(\nu-\varkappa+1)},
\end{equation}
instead of (\ref{estmulti}).

\section{Basic tools}

Let $M, N$ be integers, $M<N$, and $c_M^{\,N}=(c_{M}, c_{M+1},\ldots,c_{N})$ a nonzero sequence of complex numbers.
In this section we introduce notions  of indices and essential polynomials for the sequence $c_M^{\,N}$. 
These notions were given in more general setting in the paper~\cite{Adukov98}. Here we will consider the scalar case only. The proofs of all statements of this section can be found in~\cite{Adukov98}.

Let us form the family of all  Toeplitz matrices
\begin{equation}
T_{k}(c_M^{\,N})=
\begin{pmatrix}
 c_{k}   & c_{k-1} &\ldots & c_{M}  \\
 c_{k+1} & c_{k}   &\ldots & c_{M+1}\\
 \vdots  & \vdots  & &\vdots  \\
 c_{N}   & c_{N-1} &\ldots & c_{N+M-k}
\end{pmatrix}, \ \ \ M \leq k \leq N,
\end{equation}
which can be constructed with the help of  the sequence $c_M^{\,N}$. We will used the short 
designation $T_{k}$ in place of  $T_{k}(c_M^{\,N})$ if there is not the possibility of misinterpretation.

Our nearest aim is to 
describe a structure of the  kernels ${\ker}\,T_{k}$.
It is more convenient to deal not with  vectors $Q=(q_0, q_1, \ldots, q_{k-M})^t \in {\ker}\,T_{k}$ but with their generating polynomials $Q(z)=q_{0}+q_{1}\,z+\cdots+q_{\,k-M}\,z^{k-M}$.  We will use the spaces ${\cal N}_{k}$ of the generating polynomials instead of the spaces $\ker T_k$.
The generating function $\sum_{j=M}^Nc_kz^k$ of the sequence $c_M^{\,N}$ will be denoted by $c_M^{\,N}(z)$.

Let us introduce a linear functional $\sigma$  by the formula:
$$
\sigma\{z^j\}=c_{-j}, \ \ -N\leq j\leq -M.
$$
The functional is defined on the space of rational functions of the form 
$$
Q(z)=\sum_{j=-N}^{-M}q_jz^j.
$$
Besides this algebraic definition of $\sigma$ we will use the following analytic   
definition
\begin{equation}\label{ansigma}
\sigma\{Q(z)\}=\frac{1}{2\pi i}\int_\Gamma t^{-1}c_M^{\,N}(t)Q(t)\ dt.
\end{equation}
Here $\Gamma$ is any closed contour around the point $z=0$.

Denote by ${\cal N}_{k} \ (M \leq k \leq N)$
the space of polynomials $Q(z)$ with the formal degree
$k-M$ satisfying the orthogonality conditions:
\begin{equation}  \label{e-kerR}
\sigma \bigl\{z^{-i}Q(z)\bigr\}=0,\ i=k,k+1,\ldots,N.
\end{equation}

It is easily seen that  ${\cal N}_k$ is the space of generating polynomials  of vectors in $\ker T_k$. For convenience, we put
${\cal N}_{M-1}=0$ and denote by ${\cal N}_{N+1}$ the $(N-M+2)$-dimensional space of all polynomials with the formal degree $N-M+1$. If necessary, the more detailed notation ${\cal N}_k(c_M^{\,N})$ instead of ${\cal N}_k$ is used.

Let $d_{k}$ be the dimension of the space ${\cal N}_{k}$ and
$\Delta_{k}=d_{k}-d_{k-1}\ (M\leq k \leq N+1)$. The following proposition is crucial for the further considerations.

\begin{prop}  \label{t-sker}
For any non-zero sequence $c_{M}^{\,N}$  the following inequalities
\begin{equation}
0 = \Delta_{M}\leq \Delta_{M+1} \leq
\ldots \leq \Delta_{N} \leq \Delta_{N+1} = 2 \label{e-rchain} 
\end{equation}
are fulfilled. $\ \ \Box$
\end{prop}    

It follows from the inequalities (\ref{e-rchain}) that there exist
integers $\mu_{1}\leq \mu_{2}$ such that
\begin{equation}   \label{e-DeltaR}
\begin{array}{ccccccl}
\Delta_{M}&=&\ldots&=&\Delta_{\mu_{1}}&=&0, \\
\Delta_{\mu_{1}+1}&=&\ldots&=&\Delta_{\mu_{2}}&=&1, \\
\Delta_{\mu_{2}+1}&=&\ldots&=&\Delta_{N+1}&=&2.
\end{array}
\end{equation}
If the second row in these relations is absent, we assume $\mu_1=\mu_{2}$. 

\begin{newdef} \label{indices}
The integers $\mu_1,\ \mu_{2}$ defined in~(\ref{e-DeltaR}) will be called
 the essential indices (briefly, indices) of the sequence $c_{M}^{\,N}$.
\end{newdef}

The following proposition gives  formulas for the indices. 
\begin{prop} \label{indformula}
Let $\pi={\rm rank}\,T_{\left[\frac{N+M}{2}\right]}$, where 
$\left[\frac{N+M}{2}\right]$ is the integral part of $\frac{N+M}{2}$. Then the indices $\mu_1,\ \mu_2$
are found by the formulas:
$$
\mu_1=M+\pi-1, \ \ \ \mu_2= N-\pi+1. \ \ \ \ \ \Box
$$
\end{prop}

It follows from the definition of ${\cal N}_{k+1}$  that ${\cal N}_k$ and $z{\cal N}_k$ are subspaces of ${\cal N}_{k+1}$, $M-1\leq k\leq N$. 
Let $h_{k+1}$ be the dimension of any complement  ${\cal
H}_{k+1}$ of the subspace ${\cal N}_{k} + z{\cal
N}_{k}$ in the whole space ${\cal N}_{k+1}$.

From~(\ref{e-DeltaR}) we see that
$h_{k+1}\neq 0$ iff $k=\mu_j\ (j=1,2)$,
$h_{k+1}=1$ if $\mu_1<\mu_2$, and $h_{k+1}=2$ for $\mu_1=\mu_2$.   Therefore, for $k\neq \mu_j$
\begin{equation}
{\cal N}_{k+1}={\cal N}_k+z{\cal N}_k,
\label{e-rker1}
\end{equation}
and for $k = \mu_{j}$
\begin{equation}
{\cal N}_{k+1} = \bigl({\cal N}_k+z{\cal
N}_k\bigr) \oplus {\cal H}_{k+1}.  \label{e-rker2}
\end{equation}

\begin{newdef}\label{esspoly}
Let $\mu_1=\mu_2$. 
Any polynomials $Q_1(z)$, $Q_2(z)$ that form a basis for the two-dimensional
space ${\cal N}_{\mu_1+1}$ are called the essential
polynomials of the sequence $c_{M}^{\,N}$ corresponding
to the index $\mu_1=\mu_2$.

If $\mu_1<\mu_2$, then 
any polynomial $Q_j(z)$ that is a basis for an one-dimensi\-onal
complement ${\cal H}_{\mu_j+1}$ is said to be the essential
polynomial of the sequence  corresponding
to the index $\mu_j$, $j=1,2$.

\end{newdef}

It follows from Theorem 4.1 of the work~\cite{Adukov98}  that in the scalar case
the following criterion of essentialness is fulfilled.

\begin{prop}\label{test}
Integers
$\mu_1,\mu_2$, $\mu_1+\mu_2=M+N$,  are the indices and polynomials  $Q_1(z)\in {\cal
N}_{\mu_1+1}$, $Q_2(z)\in {\cal N}_{\mu_2+1}$ are the essential
polynomials of the sequence $c_{M}^{\,N}$  iff
\begin{equation} \label{sigma0}
\sigma_0:=\sigma\{z^{-\mu_1}Q_{2,\mu_2-M+1}Q_1(z) - 
z^{-\mu_2}Q_{1,\mu_1-M+1}Q_2(z)\}\ne 0.
\end{equation}
Here $Q_{j,\mu_j-M+1}$ is the coefficient of $z^{\mu_j-M+1}$ 
in the polynomial $Q_j(z)$.    $\ \ \Box$
\end{prop}

Now we can describe the structure of the kernels of the matrices $T_k$
in terms of the indices
and the essential polynomials.

\begin{prop}   \label{t-basis}

Let $\mu_1,\mu_{2}$ be the indices and
$Q_{1}(z),Q_{2}(z)$  the essential polynomials of the sequence
$c_{M}^{\,N}$.

Then
\begin{center}
$ {\cal N}_{k}=
\left\{
\begin{array}{c}
\{0\}, \\
\{q_{1}(z) Q_{1}(z)\}, \\
\{q_{1}(z) Q_{1}(z)+q_{2}(z) Q_{2}(z)\},
\end{array} \right.
$ $
\begin{array}{c}
M \le k \le \mu_{1},  \\
\mu_{1}+1 \le k \le \mu_{2}, \\
\mu_{2}+1 \le k \le N,
\end{array}
$
 \end{center}
where $q_{j}(z)$ is an arbitrary polynomial of the formal degree
 $k-\mu_{j}-1$, $j=1 ,2$. $\ \ \Box$
\end{prop}

\section{Construction of Wiener -- Hopf factorization of scalar polynomials}\label{sec4}

In this section we propose a method for solving the problem of the polynomial factorization in terms of indices and essential polynomials of some sequence.

Let $p(z)=p_0+p_1 z+\cdots+p_{\nu} z^{\nu}$ be a polynomial of  degree $\nu>1$, 
$p_0\ne 0$.
We suppose that $p(z)\ne 0$ on the unit circle $\mathbb T$, hence $p(z)\ne 0$ on a  closed annulus
$K:=\left\{ z\in \mathbb C : r\leq |z|\leq R\right\}$ for some $0<r<1<R<\infty$. Let 
$\varkappa = {\rm ind}_{\mathbb T}\, p(z)$ be the index of $p(z)$,
i.e. the number of zeros of the polynomial inside the unit circle. We can assume that
$0<\varkappa<\nu$,  otherwise the factorizations are trivial. Put 
$n_0=\max\{\varkappa,\, \nu-\varkappa\}$.

We will find the factorization of $p(z)$ in the form
\begin{equation}\label{factp}
p(z)=p_1(z)p_2(z),
\end{equation}
where $\deg{p_1}=\varkappa$, $\deg{p_2}=\nu-\varkappa$, the polynomial $p_1(z)$ is monic, and all zeros of $p_1(z)$ (respectively $p_2(z)$) lie in the domain $\{z\in \mathbb C : |z|<r\}$ (respectively in $\{z\in \mathbb C : |z|>R\}$). This means that 
$$
p(z)={p_-(z)}\,z^\varkappa p_+(z), \  \ z\in \mathbb T,
$$
is the Wiener -- Hopf factorization of $p(t)$ normalized by the condition $p_-(\infty)=1$.
Here $p_-(z)=\frac{p_1(z)}{z^\varkappa}$, $p_+(z)=p_2(z)$.
Let 
$$
{p^{-1}(z)} = \sum_{k=-\infty}^{\infty}c_kz^k 
$$
be the Laurent series for analytic function ${p^{-1}(z)}$ in the annulus $K$. Here
\begin{equation}\label{coefL}
c_k=\frac{1}{2\pi i}\int_{|z|=\rho}{t^{-k-1}p^{-1}(t)}\,dt, \ k\in \mathbb Z,  \ r\leq \rho\leq R.
\end{equation}
Note that $c_k$ are the Fourier coefficients of the function 
${p^{-1}(t)}$, $t\in \mathbb T$.

Form the sequence $c_{-n-\varkappa}^{\ n-\varkappa}=
 \left(c_{-n-\varkappa},\ldots,c_{-\varkappa},\ldots,c_{n-\varkappa}\right)$
for the given $n\geq n_0$. The generating function 
$c_{-n-\varkappa}^{\ n-\varkappa}(z)$ of the sequence is a partial sum of 
the Laurent series of ${p^{-1}(z)}$ and in formula~(\ref{ansigma}) we can replace 
$c_{-n-\varkappa}^{\ n-\varkappa}(z)$ by $p^{-1}(z)$.

\begin{thm}\label{esspfact}
For any $n\geq n_0$ the sequence $c_{-n-\varkappa}^{\ n-\varkappa}$  is non-zero;
the integers $-\varkappa,-\varkappa$ are the indices; and $Q_1(z)=z^{n-\varkappa+1}p_1(z)$, $Q_2(z)=p_2(z)$ are essential polynomials of the sequence.
\end{thm}
\begin{pf} Let us prove the first statement of the theorem. Since
$Q_1(z)=z^{n+1}+\alpha_{n}z^{n}+\cdots+\alpha_{n-\varkappa+1}z^{n-\varkappa+1}$, then 
$$
\sigma\{z^{-n+\varkappa-1}Q_1(z)\}=c_{-\varkappa}+c_{-\varkappa+1}\alpha_{n}+\cdots+
c_{0} \alpha_{n-\varkappa+1}.
$$
Otherwise, by analytic definition~(\ref{ansigma}) of $\sigma$
  we have 
\begin{gather*}
\sigma\{z^{-n+\varkappa-1}Q_1(z)\}=\\
\frac{1}{2\pi i}\int_{\mathbb T} t^{-1}p^{-1}(t)p_1(t)\,dt=
\frac{1}{2\pi i}\int_{\mathbb T} t^{-1}p_2^{-1}(t)\,dt =p_2^{-1}(0)\ne 0.
\end{gather*}

Therefore,
$$
c_{-\varkappa}+c_{-\varkappa+1}p_{-1}+\cdots+ c_{0} p_{-\varkappa}\ne 0,
$$ 
and not all of the numbers $c_{-\varkappa}, c_{-\varkappa+1},\ldots, c_0$ are zero. 
For $n\geq n_0$ the sequence $\{c_{-\varkappa}, c_{-\varkappa+1},\ldots, c_0\}$ is a subsequence of $c_{-n-\varkappa}^{\ n-\varkappa}$. Hence, 
$c_{-n-\varkappa}^{\ n-\varkappa}$ is also non-zero. 

Now we prove that $Q_1(z), Q_2(z)$ belong to the space ${\cal N}_{-\varkappa+1}$.
A formal degree of polynomials from this space should be $n+1$. In our case we have 
$\deg Q_1(z)=n+1$ and $\deg Q_2(z)=\nu-\varkappa<n+1$. Verify orthogonality 
condition~(\ref{e-kerR}). We have
$$
\sigma\{z^{-j}Q_1(z)\}=\frac{1}{2\pi i}\int_{\mathbb T}t^{n-\varkappa-j}p_2^{-1}(t)\,dt=0
$$
for $n-\varkappa-j\geq 0$. Here we use analytic definition~(\ref{ansigma}) of 
$\sigma$ and the Wiener -- Hopf factorization~(\ref{factp}). In particular,
$$
\sigma \bigl\{z^{-j}Q_1(z)\bigr\}=0,\ j=-\varkappa+1,\ldots,n-\varkappa.
$$
Hence $Q_1(z)\in {\cal N}_{-\varkappa+1}$.

Similarly, 
$$
\sigma\{z^{-j}Q_2(z)\}=\frac{1}{2\pi i}\int_{\mathbb T}t^{-j-1}p_1^{-1}(t)\,dt=
\frac{1}{2\pi i}\int_{\mathbb T}t^{-\varkappa-j-1}p_-^{-1}(t)\,dt=0
$$
for $j=-\varkappa+1,\ldots,n-\varkappa$, and $Q_2(z)\in {\cal N}_{-\varkappa+1}$ since
$p_-(z)=z^{-\varkappa}p_1(z)$ is analytic in the domain $\{z\in\mathbb C : |z|>1\}$.

Moreover,
$$
\sigma\{z^\varkappa Q_2(z)\}=p_-^{-1}(\infty)=1,
$$
and $Q_2(0)=p_2(0)\ne 0$.
Now we apply Propositiom~\ref{test}. Let us find the  test-number $\sigma_0$
from formula~(\ref{sigma0}):
$$
\sigma_0=\sigma\{z^{\varkappa}Q_{2,n+1}Q_1(z) - 
z^{\varkappa}Q_{1,n+1}Q_2(z)\}= -\sigma\{z^{\varkappa}Q_2(z)\}=-1\ne 0.
$$
Hence by Propositiom~\ref{test}, the integers $-\varkappa,-\varkappa$ are the indices; and $Q_1(z)$, $Q_2(z)$ are the essential polynomials of the sequence
$c_{-n-\varkappa}^{\ n-\varkappa}$. $\ \ \Box$
\end{pf}

In this theorem we have proved that there exist the essential polynomials 
$Q_1(z), Q_2(z)$ of the sequence 
$c_{-n-\varkappa}^{\ n-\varkappa}$ such that
the following additional properties are fulfilled:

\begin{enumerate}[(i)]

\item
$\deg Q_1(z)=n+1$, \quad $Q_1(0)=0$, \quad $Q_{1,n+1}=1$.
\item
$\deg Q_2(z)<n+1$, \quad $Q_2(0)\ne 0$, \quad $\sigma\{z^{\varkappa}Q_2(z)\}=1$.

\end{enumerate}

Vice versa,  if $Q_1(z), Q_2(z)\in {\cal N}_{-\varkappa+1}$ and satisfy conditions 
(i)--(ii), then $\sigma_0=-1$ and $Q_1(z), Q_2(z)$ are the essential polynomials of the sequence $c_{-n-\varkappa}^{\ n-\varkappa}$.
\begin{newdef}\label{deffactesspoly}
Let $n\geq n_0$. Polynomials $Q_1(z), Q_2(z)\in {\cal N}_{-\varkappa+1}(c_{-n-\varkappa}^{\ n-\varkappa})$ satisfying  conditions (i)--(ii) will be called the factorization essential polynomials of 
$c_{-n-\varkappa}^{\ n-\varkappa}$.
\end{newdef}

\begin{thm}\label{factesspoly}
The factorization essential 
polynomials are uniquely determined by conditions (i)--(ii).
Let $n\geq n_0+1$ and $R_1(z), R_2(z)$ are any essential polynomials of the sequence 
$c_{-n-\varkappa}^{\ n-\varkappa}$. Then the factorization essential polynomials
of $c_{-n-\varkappa}^{\ n-\varkappa}$ can be found by the formulas
\begin{align}\label{fep}
Q_1(z)&=-\frac{1}{\sigma_1}\bigl(R_{2,0}R_1(z)-R_{1,0}R_2(z)\bigr),\notag\\ 
Q_2(z)&=\frac{1}{\sigma_0}\bigl(R_{2,n+1}R_1(z)-R_{1,n+1}R_2(z)\bigr).
\end{align}
Here $R_{j,0}=R_j(0)$, $R_{j,n+1}$ is the coefficient of $z^{n+1}$ 
in the polynomial $R_j(z)$, $j=1, 2$, and
$$
\sigma_1=\begin{vmatrix}
R_{1,0}&R_{2,0}\\
R_{1,n+1}&R_{2,n+1}
\end{vmatrix}\ne 0,\ \ 
\sigma_0=\sigma\{z^\varkappa\bigl(R_{2,n+1}R_1(z)-R_{1,n+1}R_2(z)\bigr)\}\ne 0.
$$ 
 
\end{thm} 
\begin{pf} Let ${Q}_1(z), { Q}_2(z)$ and ${\widetilde Q}_1(z), {\widetilde Q}_2(z)$ be any couples of the factorization essential polynomials. Since $\{{Q}_1(z), {Q}_2(z)\}$ is a basis of the space 
${\cal N}_{-\varkappa+1}$, we have
$$
\widetilde Q_1(z)=\alpha_1{Q}_1(z)+\alpha_2{Q}_2(z), \quad
\widetilde Q_2(z)=\beta_1{Q}_1(z)+\beta_2{Q}_2(z).
$$
Then $0=\widetilde Q_1(0)=\alpha_1{Q}_1(0)+\alpha_2{Q}_2(0)=
\alpha_2{Q}_2(0)$. Hence $\alpha_2=0$. It follows from the conditions
$Q_{1,n+1}=Q_{2,n+1}=1$ that $\alpha_1=1$, and $\widetilde Q_1(z)={Q}_1(z)$.
In similar manner we can prove $\widetilde Q_2(z)={Q}_2(z)$.

Let $R_1(z), R_2(z)$ be any essential polynomials of the sequence $c_{-n-\varkappa}^{\ n-\varkappa}$ for $n\geq n_0+1$. Then $\sigma_0$ is the test-number of the
essential polynomials $R_1(z)$, $R_2(z)$ and $\sigma_0\ne 0$. 

Suppose that $\sigma_1=\begin{vmatrix}
R_{1,0}&R_{2,0}\\
R_{1,n+1}&R_{2,n+1}
\end{vmatrix}=0.$ Then
$$
\lambda_1\begin{pmatrix}
R_{1,0}\\
R_{1,n+1}
\end{pmatrix}+\lambda_2\begin{pmatrix}
R_{2,0}\\
R_{2,n+1}
\end{pmatrix} = 0. 
$$
Define $Q(z)=\lambda_1R_1(z)+\lambda_2R_2(z)$. By definition, $Q(0)=0$ and 
$\deg Q\leq n$. Put $Q(z)=z\widetilde Q(z)$, where $\deg \widetilde Q(z)\leq n-1$.
Since $Q(z)\in {\cal N}_{-\varkappa+1}(c_{-n-\varkappa}^{\ n-\varkappa})$,
we have 
$\widetilde Q(z)\in {\cal N}_{-\varkappa}(c_{-n-\varkappa+1}^{\ n-\varkappa-1})$.
However, if $n\geq n_0+1$, then the sequence $c_{-n-\varkappa+1}^{\ n-\varkappa-1}$
has indices $-\varkappa, -\varkappa$, and, by Proposition~\ref{t-basis}, 
${\cal N}_{-\varkappa}(c_{-n-\varkappa+1}^{\ n-\varkappa-1})=\{0\}$.
Thus, $Q(z)\equiv 0$, and the polynomials $R_1(z), R_2(z)$ are linearly dependent. 
But it is impossible and the contradiction proves the inequality $\sigma_1\ne 0$.

The polynomials $Q_1(z)$, $Q_2(z)$, defined by (\ref{fep}), belong to
${\cal N}_{-\varkappa +1}(c_{-n-\varkappa}^{\ n-\varkappa})$ and satisfy 
conditions (i)--(ii). Hence $Q_1(z)$, $Q_2(z)$ are the factorization essential 
polynomials. 
The theorem is proved. $\ \ \Box$
\end{pf}

Now we can construct the Wiener -- Hopf factorization of a polynomial with the help of the factorization essential polynomials. Note that existence and uniqueness of the factorization essential polynomials were proved under condition $n\geq n_0$.

\begin{thm}\label{WH}
Let $n\geq n_0$ and let $Q_1(z), Q_2(z)$ be the factorization essential polynomials of the sequence $c_{-n-\varkappa}^{\ n-\varkappa}$. Then the Wiener -- Hopf factorization of the polynomial $p(z)$ can be constructed by the formula
$$
p(z)=p_-(z)z^\varkappa p_+(z),
$$
where 
\begin{equation}\label{ff}
p_-(z)=z^{-n-1}Q_1(z), \quad p_+(z)=Q_2(z). 
\end{equation}
\end{thm}
\begin{pf} In Theorem~\ref{esspfact} we prove that $Q_1(z)=z^{n-\varkappa+1}p_1(z)$, $Q_2(z)=p_2(z)$ are
the factorization essential polynomials. By Theorem~\ref{factesspoly}, the factorization essential polynomials are determined uniquely. Thus formulas~(\ref{ff}) hold. $\ \ \Box$
\end{pf}

Now we can obtain the main result of the section about an explicit formulas for the factors of the Wiener -- Hopf factorization of a polynomial $p(z)$.

\begin{thm}\label{basic}
The matrices $T_{-\varkappa}\left(c_{-n-\varkappa}^{\ n-\varkappa}\right)$  are invertible for all $n\geq n_0$.

Let $n\geq n_0+1$. Denote by $\alpha=\left(\alpha_1,\ldots, \alpha_n\right)^T$ and 
$\beta=\left(\beta_0,\ldots, \beta_n\right)^T$ the solutions of the systems
\begin{equation}\label{basicequations}
T_{-\varkappa}\left(c_{-n-\varkappa+1}^{\ n-\varkappa-1}\right)\alpha=-(c_{-n-\varkappa}^{-\varkappa-1})^T,\ \ 
T_{-\varkappa}\left(c_{-n-\varkappa}^{\ n-\varkappa}\right)\beta=e_1,
\end{equation}
respectively. Here 
$e_1=\left(1,0,\ldots,0\right)^T$.

Then $\alpha_1=\ldots=\alpha_{n-\varkappa}=0$, $\beta_0\ne 0$, 
$\beta_{\nu-\varkappa+1}=\ldots=
\beta_n=0$, and the factors from the Wiener -- Hopf factorization of $p(z)$ are found by the formulas
$$
p_-(z)= z^{-\varkappa}(\alpha_{n-\varkappa+1}+\cdots +\alpha_n z^{\varkappa-1}+z^\varkappa), \ \ 
p_+(z)=\beta_0+\beta_1 z +\cdots+\beta_{\nu-\varkappa}z^{\nu-\varkappa}.
$$
\end{thm}
\begin{pf} For the first index of the sequence
$c_{-n-\varkappa}^{\ n-\varkappa}$  by Proposition~\ref{indformula} we have
$$
\mu_1= -n-\varkappa+\pi-1,
$$
where $\pi= {\rm rank}\, T_{-\varkappa}\left(c_{-n-\varkappa}^{\ n-\varkappa}\right)$. 
By Theorem~\ref{esspfact}, we have $\mu_1=-\varkappa$ for $n\geq n_0$. Hence rank of the 
$(n+1)\times (n+1)$ matrix 
$T_{-\varkappa}\left(c_{-n-\varkappa}^{\ n-\varkappa}\right)$ is equal to $n+1$.

If $\alpha=\left(\alpha_1,\ldots, \alpha_n\right)^T$ satisfies the equation
$T_{-\varkappa}\left(c_{-n-\varkappa+1}^{\ n-\varkappa-1}\right)\alpha=-(c_{-n-\varkappa}^{-\varkappa-1})^T$, then the vector $\left(0, \alpha_1, \ldots, \alpha_n, 1\right)^T$ belongs to the
space $\ker\, T_{-\varkappa+1}\left(c_{-n-\varkappa}^{\ n-\varkappa}\right)$, i.e.
$Q_1(t)=\alpha_1 z+\cdots+\alpha_n z^n+z^{n+1}\in {\cal N}_{-\varkappa+1}$.

Similarly, $\left(\beta_0,\ldots, \beta_n, 0\right)^T \in 
\ker\, T_{-\varkappa+1}\left(c_{-n-\varkappa}^{\ n-\varkappa}\right)$, and 
$Q_2(z)=\beta_0+\beta_1 z+\ldots+\beta_n z^n\in {\cal N}_{-\varkappa+1}$.
Moreover, $\sigma\{z^\varkappa Q_2(z)\}=c_{-\varkappa}\beta_0+\cdots+c_{-n-\varkappa}\beta_n=1$.
 
Let us prove  $\beta_0\ne 0$. If we suppose $\beta_0= 0$, then the 
non-zero vector $\left(\beta_1,\ldots, \beta_n, 0\right)^T$ belongs to 
$\ker\,T_{-\varkappa}\left(c_{-n-\varkappa+1}^{\ n-\varkappa-1}\right)$.
This is impossible since the matrix  $T_{-\varkappa}\left(c_{-n-\varkappa+1}^{\ n-\varkappa-1}\right)$  is invertible for $n\geq n_0+1$. 

Therefore the polynomials $Q_1(z)$, $Q_2(z)$ satisfy conditions (i)--(ii)
and they are the factorization essential polynomials of 
$c_{-n-\varkappa}^{\ n-\varkappa}$. By Theorem~\ref{WH},
$$
Q_1(z)=z^{n+1}p_-(z), \quad  Q_2(z)= p_+(z). 
$$

Hence, 
$$
\alpha_1=\ldots=\alpha_{n-\varkappa}=0, \quad \beta_{\nu-\varkappa+1}=\ldots=
\beta_n=0.
$$

The theorem is proved. $\ \ \Box$
\end{pf}

\section{Some computational aspects of the factorization problem}

\subsection{Computation of the index}

It is well known that the index $\varkappa$ of a scalar function  is stable under small perturbations of the function. In particular,  the following statement holds.
\begin{prop}\label{stabind}
If $p(z)\ne 0$ on the unit circle $\mathbb T$ and a polynomial $\widetilde p(z)$
satisfies the inequality
$$
\|p-\widetilde p\,\|< \min_{|z|=1}|p(z)|,
$$
then $\widetilde p(z)\ne 0$ on $\mathbb T$ and ${\rm ind}_{\mathbb T}\, p(z)={\rm ind}_{\mathbb T}\, \widetilde p(z)$.
\end{prop}
The proof is carried out using standard arguments. $\ \ \Box$

There are many ways to calculate the index. Numerical experiments show that it is more convenient to use formula~(12.6) from~\cite{Gakhov}.
Let $p(e^{i\varphi})=\xi(\varphi)+i\,\eta(\varphi)$, 
$\varphi\in [0,2\pi]$, where
$\xi(\varphi)$ and 
$\eta(\varphi)$ are real continuously differentiable functions on 
$[0,2\pi]$. Then
\begin{equation}\label{compind}
\varkappa=\frac{1}{2\pi}\int_{0}^{2\pi} \frac {\xi(\varphi)\eta'(\varphi)-\xi'(\varphi)\eta(\varphi)}{\xi^2(\varphi)+\eta^2(\varphi)}\,d\varphi.
\end{equation}

The integral can be computed numerically by the Gaussian quadrature method. Since $\varkappa$ is integer, the result is rounded up to the nearest integer.

\subsection{An a priori estimate of the condition number for the factorization problem}

Theorem \ref{basic} shows that  solving of the factorization problem 
is equivalent to solving of linear systems with the invertible matrix 
$T_{-\varkappa}(c_{-n-\varkappa}^{n-\varkappa})$, $n\geq n_0$.

Here we obtain an upper bound for the condition number 
$$k(T_{-\varkappa}(c_{-n-\varkappa}^{n-\varkappa}))=
\|T_{-\varkappa}(c_{-n-\varkappa}^{n-\varkappa})\|\,\|T^{-1}_{-\varkappa}(c_{-n-\varkappa}^{n-\varkappa})\|
$$ 
in terms of the given polynomial $p(z)$. 

Recall
that $n_0=\max\{\varkappa, \nu-\varkappa\}$, $m_K=\min_{z\in K} |p(z)|$, 
$\rho=\max \{r,1/R\}$.

\begin{prop}\label{boundT}
For $n\geq n_0$
\begin{equation}\label{upboundT}
\|T_{-\varkappa}(c_{-n-\varkappa}^{n-\varkappa})\|<\frac{1}{m_K}\,\frac{1+\rho}{(1-\rho)}.
\end{equation}
\end{prop}
\begin{pf} From the structure of the matrix  it is obvious that
$$
\|T_{-\varkappa}(c_{-n-\varkappa}^{n-\varkappa})\|\leq\|c_{-n-\varkappa}^{n-\varkappa}\|=
\sum_{j=-n-\varkappa}^{-1}|c_j|+\sum^{n-\varkappa}_{j=0}|c_j|,
$$
where
$$
c_j=\frac{1}{2\pi i}\int_{|z|=\bar\rho}{t^{-j-1}p^{-1}(t)}\,dt, \ r\leq \bar\rho\leq R,
\ j\in \mathbb Z.
$$

For $j\geq 0$ we choose $\bar\rho=R$. Then
$$
|c_j|\leq\frac{1}{2\pi}\int_0^{2\pi}\frac{|dt|}{R^{j+1}|p(t)|}\leq \frac{1}{R^{j}
\min\limits_{|z|=R}|p(z)|}\leq \frac{\rho^j}{m_K}.
$$
In the same way, if $j< 0$ we have
$$
|c_j|\leq \frac{r^{|j|}}{\min\limits_{|
z|=r}|p(z)|}\leq \frac{\rho^{|j|}}{m_K}.
$$

Therefore,
$$
\|T_{-\varkappa}(c_{-n-\varkappa}^{n-\varkappa})\|\leq\|c_{-n-\varkappa}^{\,n-\varkappa}\|\leq
\frac{1}{m_K}\left(\sum^{n+\varkappa}_{j=1}\rho^j+\sum^{n-\varkappa}_{j=0}\rho^j\right)
\leq \frac{1}{m_K}\frac{1+\rho}{1-\rho}. \ \ \  \ \ \Box
$$
\end{pf}

\begin{rmk}\label{roughest}
It is clear that $|c_j|\leq\frac{1}{m_1}$, where $m_1=\min_{|z|=1}|p(z)|$. Hence 
\begin{equation}\label{eq20}
\|T_{-\varkappa}(c_{-n-\varkappa}^{n-\varkappa})\|\leq\|c_{-n-\varkappa}^{n-\varkappa}\|\leq
\frac{2n+1}{m_1}.
\end{equation}
This rough estimate can be more precise for a polynomial of small degree.

Further, we will also use the estimates
\begin{equation}\label{estck}
\|c_{-n-\varkappa}^{-\varkappa-1}\|\leq \|c_{-n-\varkappa}^{\,n-\varkappa}\|<\frac{1}{m_K}\,\frac{1+\rho}{(1-\rho)}
\end{equation}
or $\|c_{-n-\varkappa}^{-\varkappa-1}\|\leq\frac{n}{m_1}$.
\end{rmk}

To obtain bounds for $\|T^{-1}_{-\varkappa}(c_{-n-\varkappa}^{n-\varkappa})\|$ we will get 
a formula for the inverse $T^{-1}_{-\varkappa}(c_{-n-\varkappa}^{n-\varkappa})$ in terms of
the factorization $p(z)=p_1(z)p_2(z)$. 

If $p_1(z) =z^\varkappa+p_{1,\varkappa-1}z^{\varkappa-1}+\cdots+p_{1,0}$ and $p_2(z) =z^{\nu-\varkappa}+p_{2,\nu-\varkappa-1}z^{\nu-\varkappa-1}+\cdots+p_{2,0}$, then we denote 
$$
p_1^{(k)}(z) =z^k+p_{1,\varkappa-1}z^{k-1}+\cdots+p_{1,\varkappa-k}, \ \ \ 0\leq k\leq \varkappa,
$$
and
$$
p_2^{(k)}(z) =p_{2,k}z^{k}+p_{2,k-1}z^{k-1}+\cdots+p_{2,0}, \ \ \ 0\leq k\leq \nu-\varkappa.
$$

Let $T^{-1}_{-\varkappa}(c_{-n-\varkappa}^{n-\varkappa})= (b_{ij})_{i,j=0}^n$ and 
$B_j(z)=\sum_{i=0}^n b_{ij}z^i$ be the generating polynomial of the $j$-th column $B_j$ of the inverse.

\begin{thm}\label{genB}
For $n\geq \nu$ we have 
\begin{equation}\label{formulaInv} 
B_j(z)=
\begin{cases}
p_1^{(j)}(z)p_2(z),& 0\leq j \leq \varkappa,\\
z^{j-\varkappa}p(z),&\varkappa+1 \leq j \leq n-\nu+\varkappa,\\
z^{j-\varkappa}p_1(z)p_2^{(n-j)}(z), &n-\nu+\varkappa+1 \leq j \leq n.\\
\end{cases}
\end{equation}
Here  the second row is absent if $n=\nu$.
\end{thm}
\begin{pf} If $n\geq \nu$, then  $T_{-\varkappa}(c_{-n-\varkappa}^{n-\varkappa})$ is invertible, and $Q_1(z)= z^{n-\varkappa+1}p_1(z)$, $Q_2(z)=p_2(z)$ are the factorization essential polynomials of the sequence $c_{-n-\varkappa}^{n-\varkappa}$ with the test-number
$\sigma_0=-1$. Let $\mathcal B(t,s) =\sum_{i,j=0}^nb_{ij}t^is^{-j}$ be the generating function of the inverse matrix $T^{-1}_{-\varkappa}(c_{-n-\varkappa}^{n-\varkappa})= (b_{ij})_{i,j=0}^n$.  By Theorem 2.2 
from \cite{Adukov86}, we have
$$
\mathcal B(t,s)=-s^{-n-1}\frac{Q_1(t)Q_2(s)-Q_1(s)Q_2(t)}{1-ts^{-1}}.
$$
Hence
$$
(t-s)\sum_{j=0}^nB_j(t)s^{n-j}=t^{n-\varkappa+1}p_1(t)p_2(s)-s^{n-\varkappa+1}p_1(s)p_2(t).
$$

Equating the coefficients of $s^j$, we arrive to  relation~(\ref{formulaInv}). $\ \  \Box$
\end{pf}

\begin{cor}\label{normInv}
For $n\geq \nu$, we have
\begin{equation}\label{estnormInv}
\|p\|\leq \|T^{-1}_{-\varkappa}(c_{-n-\varkappa}^{n-\varkappa})\|\leq \|p_1\|\,\|p_2\|
\leq \delta_0 \|p\|.
\end{equation}
\end{cor}
\begin{pf} It follows from (\ref{formulaInv}) that
$$
\|B_j\|=
\begin{cases}
\|p_1^{(j)}p_2\|,& 0\leq j \leq \varkappa,\\
\|p\|,&\varkappa+1 \leq j \leq n-\nu+\varkappa,\\
\|p_1p_2^{(n-j)}\|, &n-\nu+\varkappa+1 \leq j \leq n.\\
\end{cases}
$$

Hence, $\|T^{-1}_{-\varkappa}(c_{-n-\varkappa}^{n-\varkappa})\|=
\max_{0\leq j \leq n}\|B_j\|\geq \|p\|$. Moreover, $\|p_1^{(j)}p_2\|\leq \|p_1^{(j)}\|\,\|p_2\| \leq \|p_1\|\,\|p_2\|$, and $\|p_1p_2^{(n-j)}\|\leq \|p_1\|\,\|p_2^{(n-j)}\| \leq \|p_1\|\,\|p_2\|$. Thus the second inequality also holds.   $\ \ \Box$
\end{pf}

Both inequality (\ref{estnormInv}) are sharp for  $p(z)=p_1(z)p_2(z)$, where 
$p_1(z)$, $p_2(z)$ are polynomials with real nonnegative coefficients. In this case
$\|p\|=p(1)=\|p_1\|\,\|p_2\|$ (see also Proposition~\ref{spcase}).

Now, using (\ref{upboundT}), (\ref{eq20}), (\ref{estnormInv}), and (\ref{estmulti1}) we get the final result
\begin{cor}\label{estk}
For $n\geq \nu$
\begin{equation}\label{k2}
k(T_{-\varkappa}(c_{-n-\varkappa}^{n-\varkappa}))\leq 
\frac{\delta_0}{m_K}\,\frac{1+\rho}{(1-\rho)} \|p \|,
\end{equation}
and
\begin{equation}\label{k1}
k(T_{-\varkappa}(c_{-n-\varkappa}^{n-\varkappa}))\leq 
\frac{(2n+1)\delta_0}{m_1}\|p \|.
\end{equation}

\end{cor}

\subsection{Computation of the Laurent coefficients of analytic functions}
To realize the factorization method proposed in Section~\ref{sec4} we must calculate 
the Laurent coefficients 
$c_{-n-\varkappa},\ldots,c_{-\varkappa},\ldots,\allowbreak c_{n-\varkappa}$, 
of the function $p^{-1}(z)$ for  $n\geq n_0=\max\{\varkappa, \nu-\varkappa\}$.

In general, the coefficients can be found only approximately.
In order to do this, we will apply the method suggested by D.A. Bini and A. B\"ottcher  (see Theorem 3.3 in \cite{BB}). For future applications we will consider more general situation than in  this work. Moreover, our proof of inequality (\ref{estdif}) 
differs from the proof in the above mentioned theorem.

Let $f(z)$ be a function that analytic in the annulus $K=\{z\in\mathbb C: r\leq |z|\leq R\}$, $0<r<1<R<\infty$. 
By $f_k$ denote the Laurent coefficients of $f(z)$:
$$
f_k=\frac{1}{2\pi i}\int_{|t|=\rho}t^{-k-1}f(t)\,dt, \ \ r\leq\rho\leq R.
$$

For $\ell, k\in\mathbb Z$, $\ell\geq 2$, define
$$
\widetilde f_k(\ell)=\frac{1}{\ell}
\sum_{j=0}^{\ell-1}\frac{f(\omega_j)}{\omega_j^k},
$$
where $\omega_j=e^{\frac{2\pi i}{\ell}j}$, $j=0,\ldots,\ell-1$, are the zeros of the polynomial $z^\ell-1$.
\begin{thm}\label{estfk}
Let $M_K=\max\limits_{z\in K}|f(z)|$, $\rho=\max\{r,1/R\}$, and $\ell$ be an even positive integer. Then
\begin{equation}\label{estdif}
|f_k-\widetilde f_k(\ell)|<\frac{2M_K}{(1-\rho^\ell)}\rho^{\ell/2}
\end{equation}
for $k=-\ell/2,\ldots,0,\ldots, \ell/2$.
\end{thm}
\begin{pf} Define
$$
I=\frac{1}{2\pi i}\int_{\partial K}\frac{f(t)t^{\ell-k-1}}{t^{\ell}-1}\,dt.
$$
Calculate the integral by residue theorem. The integrand is analytic in the annulus $K$ except simple poles at the points $\omega_0,\ldots,\omega_{\ell-1}$. Since
$$
{\rm res}_{z=\omega_j}\frac{f(z)z^{\ell-k-1}}{z^{\ell}-1}=
\frac{f(\omega_j)}{\ell \omega_j^k}, 
$$
we have $I=\frac{1}{\ell}
\sum_{j=0}^{\ell-1}\frac{f(\omega_j)}{\omega_j^k}={\widetilde f}_k(\ell)$.

Therefore,
$$
f_k-\widetilde f_k(\ell)=\frac{1}{2\pi i}\int_{|t|=R}\frac{f(t)}{t^{k+1}}\biggl[1-\frac{t^\ell }{t^\ell-1}\biggr]\,dt+\frac{1}{2\pi i}\int_{|t|=r}\frac{f(t)t^{\ell-k-1}}{t^{\ell}-1}\,dt=$$
$$
-\frac{1}{2\pi i}\int_{|t|=R}\frac{f(t)}{t^{k+1}(t^\ell-1)}\,dt+\frac{1}{2\pi i}\int_{|t|=r}\frac{f(t)t^{\ell-k-1}}{t^{\ell}-1}\,dt.
$$
Then
$$
|f_k-\widetilde f_k(\ell)|\leq\frac{1}{2\pi}\int_{|t|=R}\frac{|f(t)|}{R^{k+1}|t^\ell-1|}\,|dt|+\frac{1}{2\pi}\int_{|t|=r}\frac{|f(t)|r^{\ell-k-1}}{|t^{\ell}-1|}\,|dt|\leq 
$$
$$
\frac{M_R }{R^k\min\limits_{|t|=R}|t^\ell-1|}+\frac{M_rr^{\ell-k} }{\min\limits_{|t|=r}|t^\ell-1|}.
$$
Here $M_R=\max_{|t|=R}|f(t)|$, $M_r=\max_{|t|=r}|f(t)|$.
It is easily seen that
$$
\min\limits_{|t|=R}|t^\ell-1|=R^\ell-1, \ \ \min\limits_{|t|=r}|t^\ell-1|=1-r^\ell.
$$
Hence, 
$$
|f_k-\widetilde f_k(\ell)|\leq\frac{M_R}{R^k(R^\ell-1)}+\frac{M_rr^{\ell-k}}{(1-r^\ell)}.
$$
Now, from the definitions of $M_K$ and $\rho$, it follows that
$$
|f_k-\widetilde f_k(\ell)|\leq\frac{M_K}{(1-\rho^\ell)}\bigl[\rho^{\ell+k}+\rho^{\ell-k}\bigr].
$$
If $k=-\ell/2,\ldots, \ell/2$, then $\rho^{\ell+k}\leq \rho^{\ell/2}$ and $\rho^{\ell-k}<\rho^{\ell/2}$. Estimate~(\ref{estdif}) has been obtained.  $\ \ \Box$
\end{pf}

By the theorem, in order to compute every element of the sequence 
$f_M,f_{M+1}\allowbreak,\ldots, f_N$
with the given accuracy, we have to select an appropriate number $\ell$.

\subsection{Stability of the factors $p_1(z)$, $p_2(z)$}

Recall that we consider the coefficients of the polynomial $p(z)$ as the initial data of the factorization problem. First of all we must study the sensitivity of the factors $p_1(z)$, $p_2(z)$ with respect to variations in these data. Moreover, we compute approximately the Laurent coefficients $c_k$ of the function $f(z)=p^{-1}(z)$ by the formula
\begin{equation}\label{appcoeF}
c_k\approx\widetilde c_k=\frac{1}{\ell}
\sum_{j=0}^{\ell-1}\frac{1}{p(\omega_j)\omega_j^k}.
\end{equation}
We can not consider these variations  as a perturbation of $p(z)$.
Hence  we must study the sensitivity of the factors $p_1(z)$, $p_2(z)$ with respect to change
in the Laurent coefficients separately.

For this reason, we first study the behavior of the factorization essential polynomials 
$Q_1(z)$, $Q_2(z)$ of the sequence $c_{-n-\varkappa}^{n-\varkappa}$ under small perturbations. Recall that the sequence consists of the Laurent coefficients of $p^{-1}(z)$ and has indices $-\varkappa, -\varkappa$, where $\varkappa={\rm ind}_{\mathbb T}\, p(z)$ and $n\geq n_0$. 

We will need some modification of a well known result on the absolute error 
for the solutions of linear systems. The statement can be proved by the standard method.

\begin{lem}\label{lemma3}
Let $A$ be an invertible matrix and $x=A^{-1}b$. If $\|A-\widetilde A\|\leq \frac{q}{\|A^{-1}\|}$ for some $0<q<1$, then $\widetilde A$ is invertible and for 
$\widetilde x=\widetilde A^{-1}\widetilde b$ we have
\begin{equation}\label{abserr}
\|x-\widetilde x\| \leq \frac{\|A^{-1}\|}{1-q}\left[\|A^{-1}\|\,\|A-\widetilde A\|\,\|b\|+\|b-\widetilde b\| \right]. \ \ \ \ \ \Box
\end{equation}
\end{lem}

\begin{thm}\label{pert1}
Let $n\geq n_0$ and let $\widetilde c_{-n-\varkappa}^{\,\,n-\varkappa}$ be a sequence such that
\begin{equation}\label{estic}
\|c_{-n-\varkappa}^{\,\,n-\varkappa}-\widetilde c_{-n-\varkappa}^{\,\,n-\varkappa}\|
\leq \frac{q}{\delta_0\|p\|}
\end{equation}
for some $0<q<1$. Then
\begin{enumerate}[(i)]
\item The indices of the sequence $\widetilde c_{-n-\varkappa}^{\,\,n-\varkappa}$ are 
$-\varkappa, -\varkappa$.

\item The sequence has the factorization essential polynomials $\widetilde Q_1(z)$,
$\widetilde Q_2(z)$.

\item For $n\geq n_0+1$ the following estimates are fulfilled:
$$
\|Q_1-\widetilde Q_1\|\leq \frac{\delta_0\|p\|}{1-q}\biggl[\frac{\delta_0\|p\|}{m_K}\frac{1+\rho}{1-\rho}+1\biggr]\|c_{-n-\varkappa+1}^{\,\,n-\varkappa-1}-\widetilde c_{-n-\varkappa+1}^{\,\,n-\varkappa-1}\|,
$$
$$
\|Q_2-\widetilde Q_2\|\leq \frac{\delta_0^2\|p\|^2}{1-q}\|c_{-n-\varkappa}^{\,\,n-\varkappa}-\widetilde c_{-n-\varkappa}^{\,\,n-\varkappa}\|.
$$
\end{enumerate}
\end{thm}
\begin{pf}
$\rm(i)$. Form the matrix $T_{-\varkappa}(\widetilde c_{-n-\varkappa}^{\,\,n-\varkappa})$.
If
\begin{equation}\label{estT1}
\|T_{-\varkappa}(c_{-n-\varkappa}^{\,\,n-\varkappa})-T_{-\varkappa}(\widetilde c_{-n-\varkappa}^{\,\,n-\varkappa})\|\leq\frac{q}{\|T_{-\varkappa}^{-1}(c_{-n-\varkappa}^{\,\,n-\varkappa})\|},
\end{equation}
then $T_{-\varkappa}(\widetilde c_{-n-\varkappa}^{\,\,n-\varkappa})$ is invertible. This means that the sequence $\widetilde c_{-n-\varkappa}^{\,\,n-\varkappa}$ has the indices
$-\varkappa, -\varkappa$.  Let us find the condition that guarantees the fulfillment of inequality~(\ref{estT1}).

From inequalities~(\ref{estnormInv}) we have
$$
\frac{q}{\delta_0\|p\|}\leq \frac{q}{\|p_1\|\,\|p_2\|}\leq \frac{q}{\|T_{-\varkappa}^{-1}(c_{-n-\varkappa}^{\,\,n-\varkappa})\|}.
$$
Since $\|T_{-\varkappa}(c_{-n-\varkappa}^{\,\,n-\varkappa})-T_{-\varkappa}(\widetilde c_{-n-\varkappa}^{\,\,n-\varkappa})\|\leq\|c_{-n-\varkappa}^{\,\,n-\varkappa}-\widetilde c_{-n-\varkappa}^{\,\,n-\varkappa}\|$, then from inequality~(\ref{estic}) it follows (\ref{estT1}).
Hence the indices of $\widetilde c_{-n-\varkappa}^{\,\,n-\varkappa}$ are 
$-\varkappa, -\varkappa$.

$\rm(ii)$. Let $\widetilde \alpha=(\widetilde \alpha_1,\ldots,\widetilde \alpha_n)$,
 $\widetilde \beta=(\widetilde \beta_0,\ldots,\widetilde \beta_n)$ be the solutions
 of systems~(\ref{basicequations}) for the perturbed sequence $\widetilde c_{-n-\varkappa}^{\,\,n-\varkappa}$. Then $\widetilde Q_1(z)=\widetilde \alpha_1z+\cdots+\widetilde \alpha_nz^n+z^{n+1}$ 
and $\widetilde Q_2(z)=\widetilde \beta_0+\cdots+\widetilde \beta_nz^n$ belong to 
$\mathcal N_{-\varkappa+1}(\widetilde c_{-n-\varkappa}^{\,\,n-\varkappa})$ and satisfy the conditions of Definition~\ref{deffactesspoly}.
Thus $\widetilde Q_1(z)$, $\widetilde Q_2(z)$ are the factorization essential polynomials 
of $\widetilde c_{-n-\varkappa}^{\,\,n-\varkappa}$.

$\rm(iii)$. For $n\geq n_0+1$ the inequality $\|c_{-n-\varkappa+1}^{\,\,n-\varkappa-1}-\widetilde c_{-n-\varkappa+1}^{\,\,n-\varkappa-1}\|\leq \frac{q}{\delta_0\|p\|}$ is valid, condition~(\ref{estT1}) is fulfilled, and, 
by Lemma~\ref{lemma3}, we have
\begin{gather*}
\|Q_1-\widetilde Q_1\|=\|\alpha-\widetilde \alpha\|\leq \frac{\|T_{-\varkappa}^{-1}(c_{-n-\varkappa+1}^{\,\,n-\varkappa-1})\|}{1-q}\times\\
\left[\|T_{-\varkappa}^{-1}(c_{-n-\varkappa+1}^{\,\,n-\varkappa-1})\|\,\|T_{-\varkappa}(c_{-n-\varkappa+1}^{\,\,n-\varkappa-1})-T_{-\varkappa}(\widetilde c_{-n-\varkappa+1}^{\,\,n-\varkappa-1})\|\,\|c_{-n-\varkappa}^{-\varkappa-1}\|\right.+\\
\left.\|c_{-n-\varkappa+1}^{\,n-\varkappa-1}-\widetilde c_{-n-\varkappa+1}^{\,\,n-\varkappa-1}\|\right ].
\end{gather*}
Taking into account estimates~(\ref{estnormInv}) and (\ref{estck}) we arrive at the desired inequality for $\|Q_1-\widetilde Q_1\|$.
The estimate for $\|Q_2-\widetilde Q_2\|$ can be proved analogously.  $\ \ \Box$
\end{pf}

Now we can study the behavior of the factors $p_1$, $p_2$ under small perturbations of $p(z)$.
Let $m_1=\min_{|z|=1} |p(z)|$.
By Proposition \ref{stabind}, if $\|p-\widetilde p\,\|< m_1$, then $\widetilde p(z)\ne 0$ on 
$\mathbb T$ and ${\rm ind}_{\mathbb T}\, p(z)={\rm ind}_{\mathbb T}\, \widetilde p(z)$.
Let $\widetilde p(z)=\widetilde p_1(z)\widetilde p_2(z)$ be the factorization of 
$\widetilde p\,(z)$. By $\widetilde c_j$ we denote the Laurent coefficients of 
$\widetilde p^{\,-1}(z)$. 

Now we estimate $\|p_1-\widetilde p_1\|$,
 $\|p_2-\widetilde p_2\|$ via $\|p-\widetilde p\,\|$.

\begin{thm}\label{pert2}

Let $n\geq n_0+1$. If $\|p-\widetilde p\,\|\leq \min\left\{q m_1, \frac{q(1-q)m_1^2}{(2n+1)\delta_0\|p\|}\right\}$, then
\begin{equation}\label{accforp1}
\|p_1-\widetilde p_1\|< \frac{(2n+1)\delta_0\|p\|}{(1-q)^2m_1^2}\, 
\biggl[\frac{\delta_0\|p\|}{m_K}\frac{1+\rho}{1-\rho}+1\biggr] \|p-\widetilde p\,\|, 
\end{equation}
and
\begin{equation}\label{accforp2}
\|p_2-\widetilde p_2\|< \frac{(2n+1)\delta_0^2\|p\|^2}{(1-q)^2m_1^2}\, 
\|p-\widetilde p\,\|. 
\end{equation}

\end{thm}
\begin{pf} Let us apply Theorem~\ref{pert1}. To do this we must estimate $\|c_{-n-\varkappa}^{\,n-\varkappa}-\widetilde c_{-n-\varkappa}^{\,\,n-\varkappa}\|$ via $\|p-\widetilde p\,\|$. 
We have
$$
\|c_{-n-\varkappa}^{\,n-\varkappa}-\widetilde c_{-n-\varkappa}^{\,\,n-\varkappa}\|=\sum_{j=-n-\varkappa}^{n-\varkappa}|c_j-\widetilde c_j |,
$$
and
$$
|c_j-\widetilde c_j |\leq \frac{1}{2\pi}\int_{|t|=1}|p^{-1}(t)-\widetilde p^{\,-1}(t)|\,|dt|
\leq \|p^{-1}-\widetilde p^{\,-1}\|_C.
$$
Hence,
$$
\|c_{-n-\varkappa}^{\,n-\varkappa}-c_{-n-\varkappa}^{\,\,n-\varkappa}\| \leq (2n+1)\|p^{-1}-\widetilde p^{\,-1}\|_C.
$$

For any Banach algebra $\mathfrak A$ the following inequality
$$
\|a^{-1}-\widetilde a^{\,-1}\|_\mathfrak A\leq \frac{\|a^{-1}\|^2_\mathfrak A}{1-q}\,\|a-\widetilde a\|_\mathfrak A
$$
holds if $\|a-\widetilde a\|_\mathfrak A\leq \frac{q}{\|a^{-1}\|_\mathfrak A}$ for some $0<q<1$.

In our case $\mathfrak A=C(\mathbb T)$, $\|p^{-1}\|_C=\frac{1}{m_1}$, and 
$\|p-\widetilde p\|_C\leq \|p-\widetilde p\|$. 
Hence, if $\|p-\widetilde p\|\leq qm_1$, then
$\|p^{-1}-\widetilde p^{\,-1}\|_C\leq\frac{1}{(1-q)m_1^2}\|p-\widetilde p\,\|_C\leq
\frac{1}{(1-q)m_1^2}\|p-\widetilde p\,\|$ and
$$
\|c_{-n-\varkappa}^{\,n-\varkappa}-\widetilde c_{-n-\varkappa}^{\,\,n-\varkappa}\| \leq \frac{(2n+1)}{(1-q)m_1^2}\|p-\widetilde p\,\|.
$$

If $\|p-\widetilde p\|\leq \min\{q m_1, \frac{q(1-q)m_1^2}{(2n+1)\delta_0\|p\|}\}$, then condition~(\ref{estic}) of Theorem~\ref{pert1} is fulfilled.
Applying this theorem we arrive desired statement. $ \ \ \Box$
\end{pf}

Now we consider perturbations  of the  polynomials $p_1(z)$, $p_2(z)$
caused by the approximation of the sequence $c_{-n-\varkappa}^{\,n-\varkappa}$ by $\widetilde c_{-n-\varkappa}^{\,n-\varkappa}$, where $\widetilde c_k=\frac{1}{\ell}
\sum_{j=0}^{\ell-1}\frac{1}{p(\omega_j)\omega_j^k}$.
Let $\widetilde p_1(z)$, $\widetilde p_2(z)$ be polynomials that define by 
Eq.~(\ref{basicequations}) for the sequence 
$\widetilde c_{-n-\varkappa}^{\,\,n-\varkappa}$.

\begin{thm}\label{pert3}
Let $n\geq n_0+1$, $\ell$ is an even integer such that $\ell \geq 2(n+\varkappa)$, and
\begin{equation}\label{nerav}
\frac{\rho^{\ell/2}}{1-\rho^\ell}<\frac {qm_K}{(4n+2)\delta_0\|p\|}.
\end{equation}

Then
$$
\|p_1-\widetilde p_1\|<\frac{(4n-2)\delta_0\|p\|}{(1-q)m_K}\left[\frac{\delta_0(1+\rho)\|p\|}{(1-q)(1-\rho)m_K}+1\right]\frac{\rho^{\ell/2}}{1-\rho^\ell},
$$
$$
\|p_2-\widetilde p_2\|<\frac{(4n+2)\delta_0^2\|p\|^2}{(1-q)m_K}\frac{\rho^{\ell/2}}{1-\rho^\ell}.
$$
\end{thm}
\begin{pf} By formula~(\ref{estdif}), we have $|c_k-\widetilde c_k|<
\frac{2}{m_K}\frac{\rho^{\ell/2}}{1-\rho^\ell}$. Therefore,
$$
\|c_{-n-\varkappa+1}^{\,n-\varkappa-1}-\widetilde c_{-n-\varkappa+1}^{\,\,n-\varkappa-1}\|<
\frac{4n-2}{m_K}\frac{\rho^{\ell/2}}{1-\rho^\ell}, \ \ 
\|c_{-n-\varkappa}^{\,n-\varkappa}-\widetilde c_{-n-\varkappa}^{\,\,n-\varkappa}\|<
\frac{4n+2}{m_K}\frac{\rho^{\ell/2}}{1-\rho^\ell}.
$$

If condition~(\ref{nerav}) is fulfilled, then 
$\|c_{-n-\varkappa}^{\,n-\varkappa}-\widetilde c_{-n-\varkappa}^{\,\,n-\varkappa}\|<\frac{q}{\delta_0\|p\|}$ and in order to obtain the required estimates it is sufficient to apply Theorem~\ref{pert1}.  $\ \ \Box$
\end{pf}

From this theorem it is easy to obtain
\begin{cor}
Let $\varepsilon<\frac{q}{1-q}$ and 
$$
\alpha=\varepsilon\frac{(1-q)m_K}{\delta_0\|p\|}\min\left\{(4n-2)\left(1+\frac{\delta_0\|p\|(1+\rho)}{m_K(1-\rho)}\right), (4n+2)\delta_0\|p\|\right\}. 
$$
If $\ell$ is an even  integer such that
\begin{equation}\label{esl}
\ell>2\max\left\{n+\varkappa,\frac{\log\left(\sqrt{1+\frac{1}{4\alpha^2}}+\frac{1}{2\alpha}\right)}{|\log \rho|} \right\},
\end{equation}
then $\|p_1-\widetilde p_1\|<\varepsilon$, $\|p_2-\widetilde p_2\|<\varepsilon$.
$\ \ \Box$

\end{cor}

\section{Algorithm and numerical examples}
The above results can be summarized in the form of the following algorithm. 

\vspace{1em}
\hrule
\vspace{0.2em}
\hrule
\vspace{0.2em}
\noindent{\bf Algorithm.} Wiener -- Hopf  factorization of a scalar polynomial
\vspace{0.2em}
\hrule
\vspace{0.2em}
\hrule
\vspace{0.3em}
\noindent INPUT. The coefficients of the polynomial $p(z)$, the parameter $\rho$ of the annulus $K$,  
$m_1=\min_{|z|=1} |p(z)|$, $m_K=\min_{z\in K}|p(z)|$, the given accuracy $\Delta$  for the coefficients of $p(z)$: $\|p-\widetilde p\|<\Delta$. 
\vspace{0.2em}
\hrule

\vspace{0.3em}
\noindent COMPUTATION.
\vspace{0.3em}
\begin{enumerate}
\item Compute the index $\varkappa$ of $p(z)$ by formula~(\ref{compind}). The result is rounded up to the nearest integer. 
\item Compute $\|p\|$, choose $n>\nu=\deg p$. For the sake of simplicity, put $q=1/2$.
\item Choose $\delta_0=1$, or $\delta_0=\varkappa+1$, or $\delta_0=\delta^\nu\sqrt{(\varkappa+1)(\nu-\varkappa+1)}$. Here $\delta=e^{2G/\pi}$ is Boyd's constant.
\item  Find accuracy $\varepsilon_1$, $\varepsilon_2$ for $p_1$, $p_2$ by formulas~(\ref{accforp1}), (\ref{accforp2}). Compute the theoretically guaranteed accuracy 
${\varepsilon}:=10^{-d}<\max\{\varepsilon_1,\varepsilon_2\}$.
\item Estimate the condition number $k\leq 10^{\tilde d}$ by formula~(\ref{k2}) or (\ref{k1}). Put ${\widetilde\varepsilon}:=10^{-d-\tilde d}$.
\item Find an even integer $\ell$ satisfying inequality~(\ref{esl}), where
$\varepsilon:=\widetilde\varepsilon$.
\item Form the sequence $\widetilde c_{-n-\varkappa}^{\,\,n-\varkappa}$ by 
formula~(\ref{appcoeF}).
\item Form the Toeplitz matrix $T_{-\varkappa+1}(\widetilde c_{-n-\varkappa}^{\,\,n-\varkappa})$ and find a basis $\{R_1$, $R_2\}$  of its kernel.
The last operation can be done with the help of SVD.
\item Find the factorization essential polynomials $Q_1(z)=\alpha_1z+\cdots+\alpha_nz^n+z^{n+1}$, $Q_2(z)=\beta_0+\beta_1z+\cdots+\beta_nz^n$ by~(\ref{fep}). 
\item Verify that the absolute values of the coefficients 
$\alpha_1,\ldots,\alpha_{n-\varkappa}$, $\beta_{\nu-\varkappa+1},\allowbreak\ldots,\beta_n$ are less than $\varepsilon$ and delete
these coefficients (see Theorem~\ref{basic}).
\item $\tilde p_1(z):=z^{-n+\varkappa-1}Q_1(z)$, $\tilde p_2(z):=Q_2(z)$.
\item end
\end{enumerate}
\vspace{0.2em}
\hrule
\vspace{0.2em}
\noindent OUTPUT. The coefficients $\tilde p_k^1$, $\tilde p_k^2$ of the factors $\tilde p_1(z)$, $\tilde p_2(z)$
with the guaranteed accuracy $\varepsilon$. 
\vspace{0.2em}
\hrule
\vspace{0.2em}
\hrule
\vspace{2em}

In the following examples we use the Maple computer algebra system. All computations were performed on a desktop.

The polynomial $p(z)$ in Example~\ref{excase1} satisfies the conditions of 
Proposition~\ref{spcase1} and its Wiener -- Hopf factorization is actually the spectral factorization.
\begin{ex}\label{excase1} Let $p(z)=(z+1/2)(z+1/3)\cdots(z+1/12)(z+2)(z+3)\cdots(z+12)$.
Taking into account the values of the coefficients of $p(z)$, we choose the precision 
$\tt Digits:=20$. Assume that the accuracy of the input data $\Delta$ is equal to $10^{-15}$.
We may take $\rho:=0.51$. Then $m_1=3.326340 \times 10^6$, $m_K=30.448076$.
\end{ex}
We have $\nu=22$, $\varkappa=11$, $\|p\|=20237817600$. Put $n=\nu+1=23$. 
By Proposition~\ref{spcase1}, $\delta_0=1$. The computation of the theoretically 
guaranteed accuracy $\varepsilon$ gives the following result 
$\varepsilon=0.695883\times 10^{-5}$. By formula~(\ref{k1}), we obtain 
the following estimate 
$k(T_{-\varkappa}(c_{-n-\varkappa}^{n-\varkappa}))\leq 2.859480\times 10^5$.
It follows from this that  $\widetilde\varepsilon=10^{-22}$ and we get 
$\ell=136$.

In this example the exact output is known 
$$
p_1(z)=(z+1/2)(z+1/3)\cdots(z+1/12), \ \ p_2(z)=(z+2)(z+3)\cdots(z+12).
$$

Table 1 shows the results of computations of the factors $\tilde p_1(z)$, $\tilde p_2(z)$.
It contains coefficients $\tilde p_k^1$, $\tilde p_k^2$, absolute errors 
$|\tilde p_k^1-p_k^1|$, $|\tilde p_k^2-p_k^2|$ for the coefficients $p_k^1$, $p_k^2$,
and $\|\tilde p_1-p_1\|$, $\|\tilde p_2-p_2\|$. For $\tilde p_k^1$, $\tilde p_k^2$ the number of decimal places obtained accurately is shown.

\begin{table}[ht]
\centering
\caption{Coefficients $\tilde p_k^1$, $\tilde p_k^2$ }
\begin{tabular}{|c|c|c|c|c|}
\hline
$k$&$\tilde p_k^1$&$|\tilde p_k^1-p_k^1|$&$\tilde p_k^2$&$|\tilde p_k^2-p_k^2|$\\
\hline
$0$&0&2.087675e-9&479001600.00000&1.04000e-9\\

$1$&0&1.60751e-7&
1007441280.00000&1.26000e-8\\

$2$&0&0.55114e-5&
924118272.00000&3.78100e-8\\

$3$&0.00011&1.97436e-18&
489896616.00000&5.78000e-8\\

$4$&0.00145&4.98731e-17&
167310220.00000&5.46800e-8\\

$5$&0.01300&9.54140e-17&
38759930.00000&3.62700e-8\\

$6$&0.08091&1.20571e-16&
6230301.00000&2.20830e-8\\

$7$&0.34928&1.18000e-16&
696333.00000&1.81875e-8\\

$8$&1.02274&8.87000e-17&
53130.00000&1.72958e-8\\

$9$&1.92925&4.76000e-17&
2640.00000&1.24386e-8\\

$10$&2.10321&1.40000e-17&
77.00000&2.80480e-9\\

$11$&1.00000&0&
0.999999&9.23655e-9\\
$\|\tilde p_1-p_1\|$&&0.56743e-5&&\\
$\|\tilde p_2-p_2\|$&&&&2.82246e-7\\
\hline
\end{tabular}
\label{tab1:coeff}
\end{table}

Thus $\|\tilde p_1-p_1\|=0.56743\times 10^{-5}<0.695884\times 10^{-5}=
\varepsilon$, and $\|\tilde p_2-p_2\|=2.82246\times 10^{-7}<0.695884\times 10^{-5}=\varepsilon$. We obtain $p_1(z)$, $p_2(z)$ with the desired accuracy.

The following example was taken from~\cite{BB}.  Since $p(z)$ has real coefficients $p_j$ and $p_{\nu-j}=p_j$,  we can use Proposition~\ref{spcase}
and the factorization of $p(z)$ is also the spectral factorization.

\begin{ex} Let $p(z)=\sum_{i=0}^{10}z^i+4z^5$, $\tt Digits:=20$, $\Delta=10^{-12}$. We may take $\rho=0.83$. Now 
$m_1=1.542464$, $m_K=0.062855$.
\end{ex}
In this example $\nu=10$, $\varkappa=5$, $\|p\|=15$, $n=\nu+1=11$. 
By Proposition~\ref{spcase}, $\delta_0=\kappa+1=6$. For the accuracy $\varepsilon$ we obtain t$\varepsilon=0.536458\times 10^{-4}$. From formula~(\ref{k1}) it follows the following estimate 
$k(T_{-\varkappa}(c_{-n-\varkappa}^{n-\varkappa}))\leq 1342.008991$.
This yields  $\widetilde\varepsilon=10^{-17}$ and $\ell=418$.

The computed coefficients of the factors $\tilde p_1(z)$, $\tilde p_2(z)$ are 
given by Table 2. We only indicate 5 decimal places here.

\begin{table}[ht]
\centering
\caption{Coefficients $\tilde p_k^1$, $\tilde p_k^2$}
\begin{tabular}{|c|c|c|c|c|c|c|}
\hline
$k$&0&1&2&3&4&5\\
\hline
$\tilde p_k^1$&0.23193&0.20715&
0.17674&0.14253&0.10685&1.00000\\
\hline
$\tilde p_k^2$&4.31154&0.46071&
0.61452&0.76203&0.89314&1.00000\\
\hline
\end{tabular}
\label{tab2:coeff}
\end{table}

In order to verify the computation correctness of $p_2(z)$, we can use the following relation between the factors $p_1(z)$ and $p_2(z)$ in the spectral factorization: 
$p_2(z)=z^\varkappa p_1(1/z)/p_1(0)$. For our example we have 
$\|\tilde p_2-z^\varkappa \tilde p_1(1/z)/\tilde p_1(0)\|=5.78\times 10^{-18}$.
Moreover, the residual error is $\|\tilde p_1 \tilde p_2- p\|=8.1\times 10^{-18}$.

In the next example the random polynomial $p(z)$ was generated with the help of package {\tt Random Tools}.
\begin{ex}
Let 
\begin{gather*}
p=z^{11}-\frac{17}{30}z^{10}+\frac{13}{10}z^9+\left(\frac{223}{60}+\frac{848}{135}i\right)z^8+
\left(-\frac{28}{15}+\frac{514}{135}i\right)z^7+\\
\left(-\frac{43}{60}+\frac{106}{135}i\right)z^6+\left(\frac{43}{60}+\frac{764}{135}i\right)z^5+\left(-\frac{31}{6}+\frac{68}{135}i\right)z^4+
\left(\frac{7}{3}-\frac{2}{3}i\right)z^3+\\
\left(-1+\frac{814}{135}i\right)z^2+\left(\frac{39}{10}+\frac{58}{15}i\right)z+
\left(-\frac{61}{60}+\frac{16}{9}i\right),
\end{gather*}
 $\tt Digits:=20$, $\Delta=10^{-18}$. 
Calculations show that $\rho=0.943396$,
$m_1=2.293009$, $m_K=0.241435$.
\end{ex}

For the polynomial we have $\nu=11$, $\varkappa=3$, $n=\nu+1=12$, $\|p\|=42.442968$. 
Since $p(z)$ is a polynomial of general type, $\delta_0=3663.225630$ is found by the formula $\delta_0=\delta^\nu\sqrt{(\varkappa+1)(\nu-\varkappa+1)}$. In fact, the use of this value of $\delta_0$ in the estimate $\|p_1\|\,\|p_2\|\leq \delta_0\|p \|$ makes it approximately 1,000 times worse. By this reasons, we are forced to use more accurate input data. Then the guaranteed accuracy in the output is $\varepsilon=0.653797\times 10^{-4}$. 

From formula~(\ref{k1}) we obtain the estimate 
$$
k(T_{-\varkappa}(c_{-n-\varkappa}^{n-\varkappa}))\leq 1.695132\times 10^6.
$$
Hence $\widetilde\varepsilon=10^{-26}$ and $\ell=1994$.

The computed coefficients of the factors $\tilde p_1(z)$, $\tilde p_2(z)$ are 
given by Table 3.

\begin{table}[ht]\label{tab1}
\centering
\caption{Coefficients $\tilde p_k^1$, $\tilde p_k^2$ }
\begin{tabular}{|c|c|c|}
\hline
$k$&$\tilde p_k^1$&$\tilde p_k^2$\\
\hline
$0$&$-0.099841 - 0.150475 i$&$-5.090491 - 10.133912 i$\\

$1$&$-0.236722 + 0.118527 i$&$-14.129949 + 0.552043 i$\\

$2$&$-0.385402 - 0.732498 i$&$-4.543939 + 4.838437 i$\\

$3$&$0.9999999$&$-7.958489 + 1.840704 i$\\

$4$&&$-5.515909 + 9.645327 i$\\

$5$&&$4.196252 + 7.320240 i$\\

$6$&&$0.930308 + 0.031004 i$\\

$7$&&$-0.181264 + 0.732498 i$\\

$8$&&$1.000000$\\
\hline
\end{tabular}
\label{tab3:coeff}
\end{table}
The residual error is $\|\tilde p_1 \tilde p_2- p\|=2.638787\times 10^{-17}$.

Let $\hat p_1(z)$, $\hat p_2(z)$ be the factorization factors of $p(z)$ obtained by the naive
method (via the roots of $p(z)$).Then 
$$\|\tilde p_1 -\hat p_1\|=1.3\times 10^{-10}, \ 
\|\tilde p_2 -\hat p_2\|=2.053051\times 10^{-8}.
$$

\section {Conclusion} 
We have considered the Wiener -- Hopf factorization problem for scalar polynomials with a numerical point of view. An algorithm has been given that is closed to {\bf algorithm 3} of 
D.A. Bini and A. B\"otcher~\cite{BB}. However, in contrast to {\bf algorithm 3}  our method permits to find coefficients of  both factors $p_1(z)$, $p_2(z)$ simultaneously. Moreover, 
effective estimates of $\|p_1-\tilde p_1\|$, $\|p_2-\tilde p_2\|$ via $\|p-\tilde p\|$ are obtained. These estimates allow to find $p_1(z)$, $p_2(z)$ with the guaranteed accuracy 
depending from the accuracy of the input data. We have illustrated this by examples.
\section*{References}


\end{document}